\documentclass[preprint,12pt]{elsart}

 \usepackage{graphicx}

\usepackage{amssymb}

\journal{Discrete Mathematics}

\begin{document}

\begin{frontmatter}



\title{$H$-product and $H$-threshold graphs.}


\author{Pavel Skums}
\ead{skumsp@gmail.com}

\address{Belarus State University, Mechanics and Mathematics Faculty, Independence av., 4, Minsk, Belarus
(currently at Molecular Epidemiology and Bioinformatics
Laboratory, Centers for Disease Control and Prevention, 1600
Clifton Road NE, Atlanta, GA, 30333)}

\begin{abstract}
This paper is the continuation of the research of the author and
his colleagues of the {\it canonical} decomposition of graphs. The
idea of the canonical decomposition is to define the binary
operation on the set of graphs and to represent the graph under
study as a product of prime elements with respect to this
operation. We consider the graph together with the arbitrary
partition of its vertex set into $n$ subsets ($n$-partitioned
graph). On the set of $n$-partitioned graphs distinguished up to
isomorphism we consider the binary algebraic operation $\circ_H$
($H$-product of graphs), determined by the digraph $H$. It is
proved, that every operation $\circ_H$ defines the unique
factorization as a product of prime factors. We define
$H$-threshold graphs as graphs, which could be represented as the
product $\circ_{H}$ of one-vertex factors, and the threshold-width
of the graph $G$ as the minimum size of $H$ such, that $G$ is
$H$-threshold. $H$-threshold graphs generalize the classes of
threshold graphs and difference graphs and extend their
properties. We show, that the threshold-width is defined for all
graphs, and give the characterization of graphs with fixed
threshold-width. We study in detail the graphs with
threshold-widths 1 and 2.

\end{abstract}


\end{frontmatter}

\section{Introduction}
\label{intro}

The decomposition methods are widely and fruitfully used in
different areas of combinatorics and graph theory. This paper is
the continuation of the previous research of the author and his
colleagues of the {\it canonical} or {\it algebraic} decomposition
of graphs. The idea of the canonical decomposition is to define
the binary operation on the set of graphs and to represent the
graph under study as a product of prime elements with respect to
this operation.

Before formulating the idea of the canonical decomposition, let us
give some basic definitions. All graphs considered are finite,
undirected, without loops and
 multiple edges. At the same time further in this paper the loops
 (but not multiple arcs) are allowed in digraphs.
 The vertex and the edge sets of a graph $G$ are
 denoted by $V(G)$ and $E(G)$, respectively. The vertex set and the arc set of a digraph $H$ are denoted by $V(H)$ and $A(H)$. Further,
 denote by $G[X]$ the subgraph induced by the set $X\subseteq V(G)$. For the convenience
 of reading the edges of graphs will be denoted as $uv$, and the
 arcs of digraphs - as $(u,v)$. Write $u\sim v$ (resp. $u\not\sim v$) if $uv\in E(G)$ (resp.
 $uv\not\in E(G)$).

 A graph $G$ is called {\it split}
\cite{FH77}, if its vertex set could be partitioned into a clique
$A$ and a independent set $B$. The graph $G$ is {\it bipartite},
if if its vertex set could be partitioned into two independent
sets $A$ and $B$. The vertex set of the complement of bipartite
graph could be partitioned into two cliques $A$ and $B$. The
partition $(A,B)$ in all those cases is called {\it a
bipartition}.

 If $X,Y\subseteq V(G)$, we will write $X\sim Y$ ($X\not\sim Y$)
if for every $x\in X$ and $y\in Y$ $x\sim y$ ($x\not\sim y$). Let
$N_Y(x)=\{y\in Y : y\sim x\}$.

The first variant of the canonical decomposition was introduced by
R. Tyshkevich and A. Chernyak \cite{TyshChern78} (in Russian) and
described in detail in \cite{Tysh00}. Consider {\it triads} (or
{\it splitted graphs}) $T=(G,A,B)$ where $G$ is a split graph and
$(A,B)$ is some fixed partition of the set $V(G)$ into clique $A$
and independent set $B$ ({\it bipartition}). The two triads
$T_i=(G_i,A_i,B_i)$, $i=1,2$, are isomorphic, if there exists  an
isomorphism $\beta : V(G_1) \rightarrow V(G_2)$
 of the graphs $G_1$ and $G_2$ preserving the bipartition $(\beta(A_1)=A_2,\ \beta(B_1)=B_2)$
Denote the set of all triads (graphs) up to isomorphism of triads
(graphs) by $Tr$ $(Gr)$.

The triads from $Tr$ could be considered as left operators acting
on the set $Gr$, the action of the operators is defined by the
formula

\begin{equation}\label{action}
  (H,A,B)\circ G=G\cup H+\{ax: a\in A, x\in V(G)\}.
\end{equation}

 On the set $Tr$ the action (\ref{action}) induces the associative binary algebraic
 operation (the {\it multiplication} of triads):

 \begin{equation}\label{binary}
 (G_1,A_1,B_1)\circ(G_2,A_2,B_2)=
 ((G_1,A_1,B_1)\circ G_2,\ A_1\cup A_2,\ B_1\cup B_2).
 \end{equation}

 A triad $T$ is called {\it decomposable} if it can be represented as a
 product of two triads. The graph is {\it decomposable}, if it is a product of a triad and a graph. Every triad $T$ can be represented as a product

\begin{equation}\label{decomp}
 T=T_1 \circ T_2\circ \ldots \circ T_k,\ k\ge 1,
\end{equation}
 of indecomposable triads $T_i$ (the parentheses in (\ref{decomp}) could be omitted because the operation $\circ$ is associative).
Analogously, every graph $G$ can be
 represented a product

\begin{equation}\label{decompgraph}
 G=T_1\circ T_2 \circ \ldots \circ T_k \circ G_0,\ k\ge 1,
\end{equation}

of indecomposable triads $T_i$ and indecomposable graph $G_0$. The
representations (\ref{decomp}) and (\ref{decompgraph}) are called
{\it the canonical decomposition} of the triad and the graphs,
respectively.

The most important property of the canonical decomposition is the
following unique factorization theorem:

\begin{thm}\label{un1decomp}\cite{Tysh00}

The canonical decomposition of the graph is determined
 uniquely,i.e. two graphs $G$ and $H$ with canonical decompositions (\ref{decompgraph}) and $H=S_1\circ S_2 \circ
 \ldots \circ S_l \circ H_0$ are isomorphic if and only if

 \begin{itemize}
 \item [1)] $k=l;$

 \item [2)] $T_i\cong S_i,$ $i=1,...,k$;

 \item[3)] $G_0\cong H_0$.
\end{itemize}

 \end{thm}

 The unique factorization property also holds for triads:

\begin{thm}\label{un1decomptriad}

The canonical decomposition of the triad is determined
 uniquely,i.e. two triads $T$ and $S$ with canonical decompositions (\ref{decomp}) and $S=S_1\circ S_2 \circ
 \ldots \circ S_l$ are isomorphic if and only if

 \begin{itemize}
 \item [1)] $k=l;$

 \item [2)] $T_i\cong S_i,$ $i=1,...,k$;

\end{itemize}

 \end{thm}

 The unique factorization theorems makes the canonical
 decomposition a very strong and useful tool to deal with the
 problems connected with the isomorphism. In particular, using the
 canonical decomposition the complete structural characterization
 of {\it unigraphs} (graphs defined up to
isomorphism by their degree sequences) was obtained by R.
Tyshkevich in \cite{Tysh00}. The crucial point of the method of R.
Tyshkevich was the fact, that the graph $G$ is a unigraph if and
only if all graphs in its canonical decomposition are unigraphs,
which follows from the unique factorization theorem. So, to
describe the structure of unigraphs it is enough to describe all
indecomposable split and indecomposable non-split unigraphs. The
description was found in \cite{Tysh00} using the properties of the
canonical decomposition and its connections with the degree
sequences of graphs.

Another applications of the canonical decomposition are the
characterizations and/or enumerations of matroidal
\cite{TyshChernChern88}, matrogenic \cite{Tysh84}, box-threshold
\cite{ChernTysh85}, domishold \cite{ChernChern90}, pseudo-split
graphs \cite{MafPre94}\cite{SST05} (these and another examples
could be found in monographs \cite{BLS99} and \cite{MP95}). The
very recent studies of the canonical decomposition and its
applications were carried out by M. Barrus and D. West
\cite{Bar10},\cite{BarWestSubm}. Among their results the very
elegant characterization of decomposable graphs from
\cite{BarWestSubm} should be especially mentioned: the graph $G$
is indecomposable if and only if its so-called {\it $A_4$-
structure} is connected. M. Barrus also applied the canonical
decomposition to the antimagic labelings of graphs \cite{Bar10}.

The success of the canonical decomposition stimulated author and
his colleagues to consider the following problem: how to
generalize canonical decomposition keeping all its advantages? The
most natural way to do it is to consider all triads $T=(G,A,B)$,
where $G$ is an {\bf arbitrary} graph and $(A,B)$ is some {\bf
arbitrary} partition of its vertex set. The multiplication
operations remain the same, as in the case of the canonical
decomposition. In this case the representations \ref{decomp} and
\ref{decompgraph} are called an {\it operator decomposition} of
triad and graph, respectively (the name came from the observation,
that the set of triads acts like the semigroup of operators on the
set of graphs). The operator decomposition was firstly considered
in \cite{TS01} (in Russian) and studied in detail in \cite{SST10}.

 It appears,
that in general the unique factorization theorem does not hold for
graphs, but holds for triads (up to permutations of staying
together commutative multipliers) \cite{SST10}. It is still a very
powerful property, which was confirmed by the applications of the
operator decomposition to the one of the most old and famous open
problems in graph theory -- the reconstruction conjecture.

Before formulating that results, let us introduce some notions. A
pair of graph classes $(P,Q)$ is called closed hereditary, if they
are hereditary,  $P$ is closed with respect to the operation of
join and $Q$ is closed with respect to the operation of disjoint
union. Graph $G$ is $(P,Q)$-{\it split}, if there exists a
partition $V(G)=A\cup B$ such, that $G[A]\in P$ and $G[B]\in Q$.
The set $M\subseteq V(G)$ is called a {\it homogeneous set}, if
every vertex $v\in V(G)\setminus M$ is adjacent either to all
vertices of $M$ or to none of them. Denote the sets of vertices of
the first and the second type by $A(M)$ and $B(M)$, respectively.
The main result of \cite{SST10} is the following. {\it Suppose
that the graph $G$ have a homogeneous set $M$ such that for some
closed hereditary pair of classes $(P,Q)$ $G[A]\in P$, $G[B]\in Q$
and $G[M]$ is not $(P,Q)$-split. Then $G$ is reconstructible.}
Note, that the property of the closed hereditariness of a pair
$(P,Q)$ is not very restrictive (there are many well-known graph
classes, which form such a pair), and so the reconstruction result
is rather general. Another applications of the unique
factorization theorem for triads in this area includes proof of
the reconstruction conjecture for $P_4$-disconnected and
$P_4$-tidy graphs \cite{ST09}.

The machinery behind the reconstruction results above is based on
the unique factorization theorem for the operator decomposition of
triads.

The further development of the theory of decomposition and its
applications requires further generalization. The natural next
step is the consideration of an arbitrary algebraic operation and
turning the set of graphs into semigroup with respect to this
operation. In this paper we study such operations.

\begin{figure}[h]\label{figH0}
\begin{center}
\includegraphics*[width=40mm] {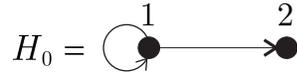}
\caption{The digraph $H_0$}%
\end{center}
\label{figHThr1}
\end{figure}

Consider the graph together with some arbitrary partition of its
vertex set into $n$ subsets. Let us call this object $n$-{\it
partitioned graph}. The isomorphism of $n$-prtitioned graphs is
naturally defined as the isomorphism of corresponding graphs
preserving the partitions. On the set of all $n$-partitioned
graphs distinguished up to an isomorphism define the binary
algebraic operation $\circ_H$ ($H$-{\it product} of graphs)
determined by the digraph $H$ with $V(H)=\{1,...,n\}$. For the two
$n$-partitioned graphs $T=(G,A_1,...,A_n)$ and $S=(F,
B_1,...,B_n)$ ($V(G)\cap V(F) = \emptyset$) their product
$S=T\circ_H S$ is the $n$-partitioned graph $(R,A_1\cup
B_1,...,A_n\cup B_n)$, where $A_i$ and $B_j$ are completely
adjacent in $F$, if $(i,j)$ is an arc of $H$, and completely
nonadjacent, otherwise. The representation of the $n$-partitioned
graph as an $H$-product of prime factors is called its {\it
$H$-decomposition}. Within this approach the operator
decomposition is $H_0$-decomposition, where the digraph $H_0$ is
shown on the figure \ref{figH0}.

The algebraic properties of the operation $\circ_H$ for 2-vertex
digraphs $H$ were studied before. The fact, that for every $H$
with $|H|=2$ the operation $\circ_H$ defines the unique
factorization of 2-partitioned graph up to the permutation of
staying together commutative multipliers, follows from the results
of \cite{LMR07}. Independently, the same fact for $H^*$ with
$A(H^*)=\{(1,2)\}$ was proved in \cite{ST07}. Moreover, in
\cite{ST07} the multiplication $\circ_{H^*}$ of a bipartite graph
with the fixed bipartition and a graph was considered (analogously
with the multiplication of a splitted graph and a graph above),
and it was proved, that in this case the unique factorization
property also holds for the decomposition of graphs, with the
exception of the simple and well-described graph family. This
unique factorization theorem was used to prove, that for the
graphs decomposable with respect to $\circ_{H^*}$ the
reconstruction conjecture is true. The last result is naturally
related to the old and well-known open problem: to prove the
reconstruction conjecture for bipartite graphs.

In fact, this kind of operations was already introduce in the
theories of clique-width \cite{CourOl00} and NLC-width
\cite{Wan94}. This two notions are similar and in some sense
equivalent, so let us quote the definition of NLC-width and the
corresponding decomposition. For a given integer $k$ consider the
set of all labeled graphs $(G,l)$, where $l$ is a mapping $l: V(G)
\rightarrow \{1,...,k\}$. The class $NLC_k$ is recursively defined
as follows \cite{Wan94}:

\begin{itemize}

\item[1)] the one-vertex labeled graphs $(K_1,l)$ belongs to
$NLC_k$;

\item[2)] if $(G_1,l_1)$, $(G_2,l_2)\in NLC_k$ ($V(G_1)\cap V(G_2)
= \emptyset$) and $S$ is some binary relation on the set
$\{1,...,k\}$, then the following labeled graph $(H,p)$ belongs to
$NLC_k$:

\begin{equation}
    V(H)=V(G_1)\cup V(G_2),
\end{equation}

\begin{equation}
    E(H)=E(G_1)\cup E(G_2)\cup \{uv : (l_1(u),l_2(v))\in S\}
\end{equation}

\begin{equation}
    p(u) = l_i(u), u\in V(G_i), i=1,2.
\end{equation}

\item[3)] If $(G,l)\in NLC_k$ and $\alpha:\{1,...,k\}\rightarrow
\{1,...,k\}$ is a function, then $(G,\alpha l)\in NLC_k$ (here
$\alpha l$ is the composition of functions).

\end{itemize}

NLC-width of a graph $G$ is the minimal $k$ such, that $G\in
NLC_k$.

Clearly, the operation in 2) is exactly the operation $\circ_H$.
But it was introduced with completely different purposes, and its
algebraic properties in general case have not been studied before.
We consider the decomposition idea from the different point of
view - as the study of binary algebraic operation. Since we want
to obtain the decomposition tool useful for the problems connected
with isomorphism (especially for the reconstrution conjecture),
the main questions, which we are interested in, is the existence
of the unique factorization property.

We also would like to note, that $H$-decomposition is related to
another well-known graph-theoretical notion --  the idea of
$M$-{\it partitions} introduced by T. Feder, P. Hell, S. Klein and
R. Motwani in \cite{FedHellKleMot03}. Suppose that $M$ is the
$n\times n$ symmetric matrix with the elements from the set
$\{0,1,*\}$. An $M$-partition of the graph $G$ is a partition
$V(G)=A_1\cup...\cup A_n$ such that each $A_i$ is either a clique
(if $M_{i,i}=1$), or independent set (if $M_{i,i}=0$), or an
arbitrary set (if $M_{i,i}=*$); and $A_i$ and $A_j$ are either
completely adjacent (if $M_{i,j}=1$), or completely nonadjacent
(if $M_{i,j}=0$), or can have arbitrary set of edges between them
(if $M_{i,j}=*$). The matrix $M$ could be considered as an
adjacency matrix of a {\it trigraph} \cite{FedHellTuck06}, which
consists of the set of $n$ vertices $\{v_1,...,v_n\}$, any two
vertices $v_i$, $v_j$ are connected either by a {\it non-edge} (if
$M_{i,j} = 0$), or {\it weak edge} (if $M_{i,j} = *$), or {\it
strong edge} (if $M_{i,j} = 1$). In this terms our decomposable
graphs are $M$-partitionable graphs, where $M_{i,i}=*$ for all $i$
and the graph formed by strong edges and non-edges of the trigraph
defined by $M$ is complete bipartite with the parts of equal size
(or, in other terms, our decomposable graphs are the graphs
admitting homomorphism to trigraphs with the above-mentioned
properties).

This paper consists of 3 parts. In the first part we define the
$H$-product. We show, that for every digraph $H$ the unique
factorization property of $H$-product of $n$-partitioned graphs
holds. Namely, for every digraph $H$ every $n$-partitioned graph
has the unique $H$-decomposition up to the permutation of staying
together commutative factors.

In the second part we define and study $H$-threshold graphs and
the corresponding dimension of graphs -- the objects based on the
binary algebraic operations defined in the first part. The idea
came both from the well-known notion of threshold graph
\cite{ChHam77} and from the theory of clique-width and NLC-width.

Threshold graphs is the important and well-studied graph class
with many interesting properties and applications. There is a
number of different equivalent definitions of threshold graph. The
most important and well-known of them are summarized in the
following theorem (those and another characterizations, properties
and applications of threshold graphs could be found in the
monograph \cite{MP95})

\begin{thm}\cite{MP95}\label{threshold}

{\it The following definitions of the threshold graph $G$ are
equivalent:}

{\it a) There exist nonnegative weights $(\alpha_v : v\in V(G))$
and a threshold $t$ such, that $U\subseteq V(G)$ is an independent
set if and only if $\sum_{v\in U}\alpha_v \leq t$.}

{\it b) There exist nonnegative weights $(\beta_v : v\in V(G))$
and a threshold $s$ such, that $uv\in E(G)$ if and only if
$\beta_u + \beta_v \geq s$.}

{\it c) For every $u,v \in V(G)$ either $N(u)\subseteq N(v)\cup
\{v\}$ or $N(v)\subseteq N(u)\cup \{u\}$}

{\it d) $G$ is split with a bipartition $(A,B)$, and the sets
$\{N_A(b) : b\in B\}$ and $\{N_B(a) : a\in A\}$ are ordered by
inclusion.}

{\it e) $G$ is $(2K_2, C_4, P_4)$-free.}

{\it f) All factors in the canonical decomposition of $G$ are
one-vertex.}

{\it g) $G$ is split with a bipartition $(A,B)$, and all factors
in the canonical decomposition of the triad $(G,A,B)$ are
one-vertex}

\end{thm}

So, according to f) and g) threshold graphs are the graphs with
the simplest canonical decompositions. The question is: what are
the simplest graphs defined by the $H$-decomposition? Following
this idea, we define $H$-threshold graphs as graphs, which could
be represented as the product $\circ_{H}$ of one-vertex factors.
We show, that every graph is $H$-threshold for some digraph $H$.
So, it is natural to look for such representation with the digraph
$H$, which is as small as possible. We define threshold-width of
the graph $G$ as a minimum size of a digraph $H$ such, that $G$ is
$H$-threshold.

The idea of threshold-width is naturally agreed with the theories
of NLC-width and cliquewidth. In particular,  the class
$lin-NLC_k$\cite{GurWan05} is defined as the set of graphs which
could be constructed by the sequence of the operations 2), 3),
where at least one multiplier is one-vertex, and linear NLC-width
of the graph $G$ is the minimal $k$ such that $G\in lin-NLC_k$. In
the case of threshold-width, the operation is fixed.

Another important graph class of graphs related to the threshold
graphs, is the class of difference graphs \cite{HPS90} (another
name is bipartite chain graphs \cite{Yan82}).

\begin{thm}\cite{HPS90}\label{differgrdef}{\it  The following definitions of the difference graph $G$ are
equivalent:}

{\it a) There exist real weights $(\beta_v : v\in V(G))$ and a
threshold $s$ such, that $|\beta_v|\leq t$, $v\in V(G)$ and $uv\in
E(G)$ if and only if $|\beta_u - \beta_v| \geq s$.}

{\it b) $G$ is bipartite with a bipartition $(A,B)$, and the sets
$\{N_A(b) : b\in B\}$ and $\{N_B(a) : a\in A\}$ are ordered by
inclusion.}
\end{thm}

In \cite{HPS90} authors emphasize, that properties of difference
graphs are very similar to properties of threshold graphs. We
show, that it is not the coincidence, because difference graphs
are $H$-threshold for the particular $H$. So, graphs with fixed
threshold-width are direct generalizations of both threshold and
difference graphs, and we show, that they extend another
properties of those classes. In particular, we show, that graphs
with fixed threshold-width are also characterized in terms of
vertex partitions into cliques and independent sets and the
orderings of vertex neighborhoods, though the characterization
become much more complicated. More precisely, we prove, that a
graph $G$ has threshold-width at most $k$ if and only if

a) $V(G)$ could be partitioned into $k$ cliques and independent
sets $V_1$,...,$V_k$;

b) for every $i,j=1,...,k$, $i\ne j$ the sets $\{N_{V_j}(b) : b\in
V_i\}$ are ordered by inclusion;

c) those orderings for different $i$ and $j$ are coordinated in
the following sense: we can associate with the orderings the graph
$R$ and the digraph $F$  such, that $R$ is bipartite and $F$ is
acyclic.

In the third part of the paper we consider the graphs with small
threshold-width. By the definition the only graphs with threshold
dimension 1 are complete and empty graphs. Threshold graphs have
threshold-width at most 2, but there are non-threshold graphs with
this property. We give the structural characterization and the
characterization by the finite list of forbidden induced subgraphs
for the class of graphs with threshold-width at most 2.

In particular, we show, that graph $G$ has threshold-width at most
2 if and only if $G$ or $\overline{G}$ is either threshold or
difference. It is interesting to compare this characterization
with the characterization of the graphs with small linear
$NLC$-width from \cite{Gur06NLC}: a graph $G$ has linear
$NLC$-width 1 if and only if $G$ is threshold.

\section{$H$-product of graphs}
\label{Hdecomp}

Let $H$ be a digraph with the vertex set $V(H) = \{1,...,n\}$ and
the arc set $A(H)$. The {\it $n$-partitioned graph} is a
$(n+1)$-tuple $T = (G,A_1,...,A_n)$, where $G$ is a graph and
$(A_1,...,A_n)$ is a partition of its vertex set into disjoint
subsets: $V(G)=A_1\cup...\cup A_n$, $A_i\cap A_j = \emptyset$ for
all $i\ne j$. Some of sets $A_i$ could be empty. $G$ is called
{\it the basic graph} of $T$. Denote the set of vertices and the
set of edges of $T$ by $V(T)$ and $E(T)$, respectively.

The isomorphism $f$ of $n$-partitioned graphs $T$ and $S
=(F,B_1,...,B_n)$ is an isomorphism of $G$ and $F$ such that
$f(A_i) = B_i$ for every $i=1,...,n$. Let $\Sigma_n$ be the set of
all $n$-partitioned graphs distinguished up to isomorphism.

On the set $\Sigma_n$ consider a binary algebraic operation
$\circ_H : \Sigma_n \times \Sigma_n \rightarrow \Sigma_n$ ({\it
$H$-product} of $n$-partitioned graphs) as follows:

\begin{equation}\label{operation}
    (G,A_1,...,A_n)\circ_H (F,B_1,...,B_n) = (R,A_1 \cup B_1,...,A_n
\cup B_n),
\end{equation}

where $V(R)=V(G)\cup V(F)$ (we assume without lost of generality
that $V(G)\cap V(F)=\emptyset$), $E(R)=E(G)\cup E(H) \cup \{xy:
x\in A_i, y\in B_j, (i,j)\in A(H)\}$.

For the convenience we will further sometimes denote the operation
$\circ_H$ simply by $\circ$, if it is clear, what digraph $H$ we
mean. The operation, which was introduced and studied in
\cite{SST10}, is the particular case of $\circ_H$ for a digraph
$H=H_0$ shown in the figure \ref{figHThr1}

It is easy to check, that for every digraph $H$ the operation
$\circ_H$ is associative. So, the set $\Sigma_n$ with the
operation $\circ_H$ is a semigroup.

The digraph $H$ is {\it symmetric}, if $(i,j)\in A(H)$ whenever
$(j,i)\in A(H)$. It is clear that the operation $\circ_H$ is
commutative if and only if $H$ is symmetric.

The $n$-partitioned graph $T\in \Sigma_n$ is called $H$-{\it
decomposable}, if $T=T_1\circ_H T_2$, $T_1,T_2\in \Sigma_n$, and
$H$-{\it prime}, otherwise. It is clear, that every
$n$-partitioned graph $T\in \Sigma_n$ could be represented as a
product $T=T_1\circ_H...\circ_H T_k$, $k\geq 1$, of $H$- prime
factors. Such a representation is called an $H$-{\it
decomposition} of $T$.

\begin{thm}\label{uniq}({\bf unique factorization theorem for the operation $\circ_H$}) {\it For every $n$-vertex digraph $H$ every $n$-partitioned graph
$T\in \Sigma_n$ has the unique $H$-decomposition up to the
permutation of staying together commutative factors.}

\end{thm}

\begin{pf} It is evident, that if
two $n$-partitioned graphs have the $H$-decompositions, which
differ only by some permutations of staying together commutative
multipliers, then they are isomorphic. So let us prove the inverse
proposition. It is evident for prime $n$-partitioned graphs.
Further apply the induction by the number of vertices.

Let

\begin{equation}\label{HdecompUV}
    U=T_1\circ...\circ T_k,\;W=R_1\circ...\circ R_l,
\end{equation}

$U\cong W$; $k,l \geq 2$. Let $$U=(G,X_1,...,X_n),\,\,
W=(F,Y_1,...,Y_n).$$ We may assume that $X_i\cup Y_i\ne \emptyset$
for all $i=1,...,n$.

Let $f: V(U)\rightarrow V(W)$ is the isomorphism of $U$ and $W$.
We will use the following notation. For the set $X\subseteq V(U)$
let $f(X)=\{f(x) : x\in X\}$, for the subgraph $G'$ of $G$ let
$f(G') = W[f(V(G'))]$ and for the $n$-partitioned graph
$T=(G',A_1,...,A_n)$, where $G'$ is a subgraph of $G$, let $f(T) =
(f(G'),f(A_1),...,f(A_n))$.

 Setting  $S=T_2\circ...\circ T_k = (G'',S_1,...,S_n)$,
$Q=R_2\circ...\circ R_l=(F'',Q_1,...,Q_n$), we have

\begin{equation}
    U=T_1\circ S,\,\, W=R_1\circ Q.
\end{equation}

Let $T_1 = (G',A_1,...,A_n)$, $R_1 = (F',B_1,...,B_n)$.  By the
definition of the isomorphism $f(A_i\cup S_i) = B_i\cup Q_i$.

Suppose that there exists $i\in \{1,...,n\}$ such that $f(A_i)\cap
B_i\ne \emptyset$, $f(A_i)\cap Q_i \ne \emptyset$. Then
$$f(T_1)=T'\circ T'',$$ where $$T' = (F[f(V(T_1))\cap V(R_1)],f(A_1)\cap
B_1,...,f(A_n)\cap B_n),$$ $$T'' = (F[f(V(T_1))\cap
V(Q)],f(A_1)\cap Q_1,...,f(A_n)\cap Q_n).$$ Here $V(T'),V(T'')\ne
\emptyset$ by the assumption. It contradicts the fact that $T_1$
is prime.

Analogously, the existence of $i\in \{1,...,n\}$ such that
$f^{-1}(B_i)\cap A_i\ne \emptyset$, $f^{-1}(B_i)\cap S_i\ne
\emptyset$ contradicts the fact, that $R_1$ is prime.

So, further we can assume that for every $i=1,...,n$
$f(A_i)\subseteq B_i$ or $f(A_i)\subseteq Q_i$.

Suppose that there exist $i,j\in \{1,...,n\},$ $i\ne j$ such that
$f(A_i)\subseteq B_i$ and $f(A_j)\subseteq Q_j$. Then
$f(T_1)=T'\circ T''$, where $T',T''$ are defined as above. Again
the contradiction with the indecomposibility of $T_1$ is obtained.

So, there are two possibilities:

1) For every $i=1,...,n$ $f(A_i)\subseteq B_i$. Then the facts
proved above imply, that $f(A_i)=B_i$, $f(S_i)=Q_i$ for every
$i=1,...,n$. Thus $T_1\cong R_1$, $S\cong Q$. After applying
induction assumption to the $S=T_2\circ...\circ T_k$ and
$Q=R_2\circ...\circ R_l$, we get, that $k=l$ and under the
respective ordering $R_2 \cong T_2$,...,$T_k\cong R_k$.

2) For every $i=1,...,n$ $f(A_i)\subseteq Q_i$. Then $B_i\subseteq
f(S_i)$.

Let $f(S_i)\cap Q_i = \emptyset$ for all $i=1,...,n$. Then
$f(S_i)=B_i$, $f(A_i)=Q_i$ for every $i=1,...,n$. It means, that
$S\cong R_1$, $T_1\cong Q$, and thus
$$U\cong T_1\circ R_1\cong W \cong R_1\circ T_1.$$ So, the
statement of the theorem is true.

Consider the case, when there exist $i\in \{1,...,n\}$ such that
$f(S_i)\cap Q_i\ne \emptyset$. Let $$Z=(F[f(V(S))\cap
V(Q)],f(S_1)\cap Q_1,...,f(S_n)\cap Q_n).$$ By the assumption
$V(Z)\ne \emptyset$. Then $f(S)=R_1\circ Z$, $Q = f(T_1)\circ Z$
and thus $$S\cong R_1\circ Z,\,\, Q\cong T_1\circ Z.$$ So, $T_1$
is the first factor in some $H$-decomposition of $Q$. Applying the
induction assumption to $Q$, we may assume without lost of
generality, that $T_1=R_2$ and $Z=R_3\circ...\circ R_l$. So,
$$T_2\circ...\circ T_k \cong S\cong R_1\circ R_3\circ...\circ R_l.$$

By the induction assumption applied to $S$, we have $k=l$ and
under the respective ordering $R_1 \cong T_2$, $T_3 \cong
R_3$,...,$T_k\cong R_k$.

To complete the proof, it remains to show, that $T_1$ and $R_1$
commutate. To do it, it is sufficient to prove, that for every
pair $i,j\in \{1,...,n\}$, $i\ne j$, such that $(i,j)\in A(H)$ and
$(j,i)\not\in A(H)$ one of the following four conditions hold:
either $A_i\cup A_j = \emptyset$, or $A_j\cup B_j = \emptyset$, or
$A_i\cup B_i = \emptyset$, or $B_i\cup B_j = \emptyset$.

We have $f(A_i)\sim f(S_j)$, $f(A_j)\not \sim f(S_i)$ (because
$A_i\sim S_j$, $A_j\not \sim S_i$ and $f$ is an isomorphism).

But then, since $f(A_i)\subseteq Q_i$, $f(A_j)\subseteq Q_j$,
$B_i\subseteq f(S_i)$, $B_j \subseteq f(S_j)$, we have $f(A_i)\not
\sim f(S_j)$, $f(A_j)\sim f(S_i)$.

This two facts imply, that one of the following is true:

1) $A_i\cup A_j = \emptyset$;

2) $A_i\cup S_i = \emptyset$, which implies, that $B_i =
\emptyset$;

3) $A_j\cup S_j = \emptyset$, which implies, that $B_j =
\emptyset$;

4) $S_i\cup S_j = \emptyset$, which implies, that $B_i =
\emptyset$, $B_j = \emptyset$.

The theorem is proved. \end{pf}

\section{$H$-threshold graphs and the threshold-width of graphs}
\label{HThresh}

Denote by $K^k_i$ the $k$-partitioned graph
$(K_1,\emptyset,...,\emptyset,\{v\},\emptyset,...,\emptyset)$ (the
only nonempty set of the partition is the $i$th set).

Let $H$ be a digraph on $k$ vertices. Let us call a graph $G$
$H$-{\it threshold graph}, if it is basic for the $n$-partitioned
graph of the form

\begin{equation}\label{ThrForm}
    K^k_{i_1}\circ_H...\circ_H K^k_{i_n}.
\end{equation}

 In this case for the
simplicity of the notation we will write
$G=K^k_{i_1}\circ_H...\circ_H K^k_{i_n}$ (though strictly speaking
the left part of this equality is the graph and the right part is
$k$-partitioned graph). The representation of the graph $G$ in the
form (\ref{ThrForm}) is called a {\it threshold representation} of
$G$.

To illustrate the notion of $H$-threshold graph, we show the
threshold representations of graphs $P_4$ and $C_4$ for different
2-vertex digraphs $H$ on the figure 2 (the 2-partitioned factors
$K^2_i$ are represented by ovals).

\begin{figure}[h]
\begin{center}
\includegraphics*[width=100mm] {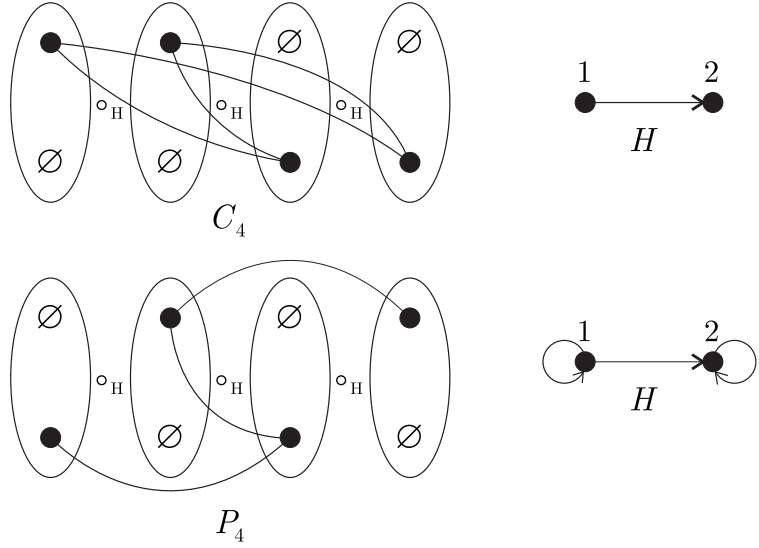}
\caption{$C_4$ and $P_4$ as $H$-threshold graphs}%
\end{center}
\label{figHThrC4P4}
\end{figure}

By Theorem \ref{threshold} threshold graphs are exactly
$H_0$-threshold graphs for the digraph $H_0$ shown in the figure
\ref{figHThr1}

\begin{prop}\label{Htrevery}{\it Every graph $G$ is $H$-threshold for some digraph $H$.}

\end{prop}

\begin{pf} Let $V(G)=\{1,...,n\}$. Define $H$ as follows:
$V(H)=\{1,...,n\}$, $(i,j)\in A(H)$ if and only if $ij\in E(G)$,
$i < j$ (i.e. $H$ is obtained from $G$ by assigning the
orientation on every edge of $G$). It is easy to see, that $G =
K^n_1\circ_H K^n_2\circ_H...\circ_H K_n^n$. \end{pf}

The digraph $H$ constructed in the proof of Proposition
\ref{Htrevery} has $|V(G)|$ vertices. But, for example, threshold
graphs are $H$-threshold for the digraph $H$ with only 2 vertices.
So, it is natural to consider the minimum order of a digraph, for
which a graph $G$ is $H$-threshold. Here we introduce the
corresponding graph parameter.

The {\it threshold-width} of a graph $G$ is the parameter
$ThrWidth(G) = min\{|H| : G \,\, is \,\,  H-threshold\}$. By the
Proposition \ref{Htrevery} every graph has the threshold -width.
It is clear, that for every graph $G$ on $n$ vertices
$ThrWidth(G)\leq n$.

\begin{prop}\label{obs1}
{\it For every graph $G$ $ThrWidth(G) = ThrWidth(\overline{G})$.}

\end{prop}

\begin{pf}
Suppose, that $G$ is $H$-threshold for a digraph $H$ with the
vertex set $V(H)=\{1,...,k\}$, i.e. $G=K_{i_1}^k\circ_H...\circ_H
K_{i_n}^k$. Let $\{v_j\} = V(K^k_{i_j})$, $j=1,...,n$. Consider
the vertices $v_p$ and $v_q$. Suppose, that $p < q$. Then $v_p\sim
v_q$ if and only if one of the following conditions hold:

1) $i_p=i_q$ and $(i_p,i_p)\in A(H)$;

2) $i_p\ne i_q$ and $(i_p,i_q)\in A(H)$.

Define $\overline{H}$ be the complement of $H$, i.e. the digraph
with the same vertex set and with the arc set $A(\overline{H}) =
\{(i,j) : (i,j)\not\in A(H)\}$. Then
$\overline{G}=K_{i_1}^k\circ_{\overline{H}}...\circ_{\overline{H}}
K_{i_n}^k$, where $\{v_j\} = V(K^k_{i_j})$, $j=1,...,n$.\end{pf}

Now we are going to give the characterization of graphs with
$ThrWidth(G)\leq k$. But firstly we need some auxiliary
definitions and lemmas.

For a digraph $H$ and $v\in V(H)$ let $N_{in}(v) = \{u\in
V(H)\setminus \{v\}: (u,v)\in A(H)\}$ and $N_{out}(v) = \{w \in
V(H)\setminus \{v\} : (v,w)\in A(H)\}$ be the {\it
in-neighborhood} and the {\it out-neighborhood} of $v$,
respectively.

Let $H$ be a digraph and let $(v_1,...,v_n)$ be the ordering of
its vertices. This ordering is called {\it acyclic ordering} or
{\it topological sort}, if all arcs of $H$ have the form
$(v_i,v_j)$, where $i < j$. A digraph is {\it acyclic}, if it does
not contain directed cycles. The following property of acyclic
graphs is well-known.

\begin{prop}\cite{Gutin}\label{acyclicordering}

{\it A digraph is acyclic if and only if there exists an acyclic
ordering of its vertices.}

\end{prop}

Let $S$ be the family of sets
$S=(\{X_1^1,X_2^1\},...,\{X_1^n,X_2^k\})$, where $X^i_j \subseteq
\{1,...,k\}\setminus \{i\}$, $i=1,...,k$, $j=1,2$ (some of sets
$X^i_j$ could be empty). Let us call $S$ a {\it digraphical
family}, if there exists a digraph $D$ on the vertex set
$V(D)=\{1,...,k\}$ such, that
$S=(\{N_{in}(1),N_{out}(1)\},...,\{N_{in}(k),N_{out}(k)\})$. $D$
is called {\it a realization} of $S$.

The evident necessary condition for the digraphicity of $S$ is
$i\in X_1^j\cup X_2^j$ whenever $j\in X_1^i\cup X_2^i$. Let us
call the family $S$ with this property {\it proper}.

Suppose that $S$ is the proper family. Define the graph $R(S)$ as
follows:  $V(R(S)) = S$, $X^i_q\sim X^j_p$ if and only if either
$i=j$, $q\ne p$ or $i\in X^j_p$, $j\in X^i_q$, $i,j=1,...,k$,
$q,p=1,2$.

\begin{lem}\label{digraphical} {\it The proper family $S$ is digraphical if and only if the graph
$R(S)$ is bipartite.}

\end{lem}

\begin{pf} Suppose that $D$ is a realization of $S$. Let

$$l(X^i_q)=\left \{
\begin{array}{ll}
1,\; if\; X^i_q = N_{out}(i) \\ 2,\; if\; X^i_q = N_{in}(i).
\end{array}
\right. $$

By the definition $l(X^i_1)\ne l(X^i_2)$, $i=1,...,k$. If $j\in
X^i_q = N_{in}(i)$, then $i\in X^j_p = N_{out}(j)$, and so
$l(X^i_q)\ne l(X^j_p)$. This $l$ is a proper 2-coloring of $R(S)$.

Inversely, let $l$ be a proper 2-coloring of $R(S)$. Define the
digraph $D$ on the vertex set $\{1,...,k\}$ as follows: $(i,j)\in
A(G)$ if and only if $i\in X^j_p$, $j\in X^i_q$, $l(X^i_q)=1$,
$l(X^j_p)=2$.

Since $l$ is a proper 2-coloring, this definition correctly
defines a digraph, and for every $i=1,...,k$ if, for example,
$l(X^i_1)=1$, $l(X^i_2)=2$, then $X^i_1 = N_{out}(i)$, $X^i_2 =
N_{in}(i)$. \end{pf}

\begin{cor}\label{observrealiz} {\it If $D_1$ and $D_2$ are two different realizations of $S$,
then $D_1$ could be obtained from $D_2$ by the reversal of all
arcs of some of its connected components.}

\end{cor}

For a sequence $\pi = (a_1,...,a_n)$ denote by $inv(\pi)$ the
sequence $(a_n,...,a_1)$.

Let

\begin{equation}\label{part}
    V(G) = V_1\cup...\cup V_k
\end{equation}

is a partition of the vertex set of the graph $G$, where each
$V_i$ is either a clique or an independent set.

We will say, that the partition (\ref{part}) satisfies the {\it
neighborhoods ordering} property, if for every $i=1,...,k$ there
exists a permutation $\psi(i) = (u_1^i,...,u_{r_i}^i)$ of the set
$V_i$ such, that for every $j\in [k]\setminus\{i\}$ the set
$\{N_{V_j}(u) : u\in V_i \}$ is ordered by inclusion and this
ordering either coincides with $\psi(i)$ or with $inv(\psi(i))$.
In other words, for every $j\in [k]\setminus\{i\}$ either

\begin{equation}\label{firsttype}
    N_{V_j}(u_1^i)\supseteq N_{V_j}(u_2^i)\supseteq...\supseteq
N_{V_j}(u_{r_i}^i)
\end{equation}

or

\begin{equation}\label{secondtype}
N_{V_j}(u_1^i)\subseteq N_{V_j}(u_2^i)\subseteq...\subseteq
N_{V_j}(u_{r_i}^i)
\end{equation}

Assume, that the permutations $\psi(i)$ are fixed. For every $i\in
[k]$ the set $[k]\setminus\{i\}$ is partitioned into two classes.
Let us for convenience denote those classes $Y^i_1$ (contains $j$
satisfying (\ref{firsttype})) and $Y^i_2 =
([k]\setminus\{i\})\setminus Y^i_1$ (contains $j$ satisfying
(\ref{secondtype}))

Let $$X^i_r = Y^i_r\setminus \{ j : V_i\sim V_j\; or\; V_i\not\sim
V_j\},\,\, r=1,2,$$ and $$S=S(V_1,...,V_k) =
(\{X_1^1,X_2^1\},...,\{X_1^k,X_2^k\}).$$

Suppose that $S$ is a digraphical family (i.e. by the Lemma
\ref{digraphical} $R(S)=R(V_1,...,V_k)$ is a bipartite graph) and
$D$ is its realization, $V(D)=[k]$. Assume without lost of
generality, that $N_{out}(i)=X_1^i$ (if it is not the case,
replace $\psi(i)$ by $inv(\psi(i))$). Note also, that by the
definition for every $i,j\in [k], i\ne j$ $D$ contains at most one
arc from the set $\{(i,j),(j,i)\}$.

Using the digraph $D$, define the digraph $F = F(V_1,...,V_k)=
F_{D}(V_1,...,V_k)$ as follows: $V(F) = V(G)$, $A(F) = A_1\cup
A_2\cup A_3^1\cup...\cup A_3^k$, where

$$ A_1 = \{(u,v) : u\in V_i, v\in V_j, uv\in E(G), (i,j)\in
A(D);i,j\in [k]\};$$

$$ A_2 = \{(v,u) : u\in V_i, v\in V_j, uv\not\in E(G), (i,j)\in
A(D); i,j\in[k]\};$$

$$ A_3^i = \left \{
\begin{array}{ll}
\{(u^i_1,u^i_2),...,(u^i_{r_i-1},u^i_{r_i})\},\; if\; X^i_1 =
N_{out}(i)\, in \, D \\ \{(u^i_{r_i},u^i_{r_{i-1}}),...,(u^i_2,
u^i_1)\},\; if\; X^i_2 = N_{out}(i) \, in \, D.
\end{array}
\right. i=1,...,k.$$

In other words, the digraph $F$ is constructed in the following
way. Firstly consider every pair $V_i,V_j$ such, that neither
$V_i\sim V_j$ nor $V_i\not\sim V_j$. Without lost of generality
suppose, that $(i,j)\in A(D)$. Consider the set $E_{i,j}$ of edges
of the complete bipartite graph with the parts $V_i$ and $V_j$. If
the edge $uv\in E_{i,j}$ belongs to $E(G)$, then orientate it in
the direction from $V_i$ to $V_j$; otherwise orientate it in the
direction from $V_j$ to $V_i$. Next turn every set $V_i$ into the
oriented path, the order of vertices of this path is defined
either by $\psi(i)$ or by $inv(\psi(i))$ (depending on what of the
sets $X^i_1$ or $X^i_2$ is the out-neighborhood of $i$ in $D$).

Now we are ready to formulate the characterization of graphs with
the threshold-width $ThrWidth(G)\leq k$.

\begin{thm}\label{HThreshcharact}

{\it Let $G$ be a graph. $ThrWidth(G)\leq k$ if and only if there
exists a partition (\ref{part}) such that}

\begin{itemize}
{\it \item[1)] it satisfies the neighbourhoods ordering property;}

\item[2)] {\it the family $S=S(V_1,...,V_k)$ is digraphical (i.e.
the graph $R(S)=R(V_1,...,V_k)$ is bipartite);}

\item[3)] {\it the digraph $F = F(V_1,...,V_k)$ is acyclic.}
\end{itemize}

\end{thm}

\begin{pf} Let us prove sufficiency first. Suppose, that $D$ is a
realization of $S$, which defines $F$. Let us expand $D$ by adding
the set of arcs $\{(i,i) : V_i \;\; is\;a\;clique\}\cup
\{(i,j),(j,i) : V_i\sim V_j\}$. Denote the obtained graph by $H$.

Let $(v_1,...,v_n)$ be an acyclic ordering of the digraph $F$. We
will show, that $G=K_{i_1}^k\circ_H...\circ_H K_{i_n}^k$, where
$V(K_{i_j}^k) = \{v_j\}$, $v_j\in V_{i_j}$.

Let $K_{i_1}^k\circ_H...\circ_H K_{i_n}^k = Z$. Consider the edge
$ab\in E(G)$. Let us show, that $ab\in E(Z)$. If $V_i$ is a clique
in $G$, then $(i,i)\in A(H)$, which implies, that $V_i$ is a
clique in $Z$. Analogously, if $V_i\sim_G V_j$, then
$(i,j),(j,i)\in A(H)$, and so by the definition of the operation
$\circ_H$ $V_i\sim_Z V_j$.

So, it remains to consider the case, when $a\in V_i$, $b\in V_j$,
$i\ne j$ and neither $V_i\sim V_j$ nor $V_i\not\sim V_j$. In this
case $i$ and $j$ are connected by an arc in $D$. Let without lost
of generality $(i,j)\in A(D)$. So $(a,b)\in A(F)$ by the
definition of $F$. Then in the acyclic ordering $a$ goes before
$b$, i.e. $a=V(K_{i_r}^k)$, $b=V(K_{i_s}^k)$, $r < s$. It together
with the fact, that $(i,j)\in A(H)$, implies that $ab\in E(Z)$.

Conversely, let $ab\in E(Z)$. Let $a=v_r$, $b=v_s$, $r < s$ (i.e.
$a$ precedes $b$ in the acyclic ordering), $a\in V_i$, $b\in V_j$.
So we know, that $V_i\not\sim V_j$ does not hold. By the
definition of the operation $\circ_H$ $(i,j)\in A(H)$. If $i=j$,
then $V_i$ is a clique, and so $ab\in V(G)$. So let further $i\ne
j$ and it is not true, that $V_i\sim V_j$. Then $(i,j)\in A(D)$
and so $a$ and $b$ are adjacent in $F$. Since $a$ precedes $b$ in
the acyclic ordering, $(a,b)\in A(F)$. So the arc $(a,b)$ is
directed from $V_i$ to $V_j$, which implies, that $ab\in E(G)$.

Now we will prove necessity. Assume, that
$G=K_{i_1}^k\circ_H...\circ_H K_{i_n}^k$, where $\{v_j\} =
V(K^k_{i_j})$. Then

\begin{equation}\label{part2}
    V(G) = V_1\cup...\cup V_k.
\end{equation}

where $V_i = \{v : \{v\} = V(K^k_{i})\}$, $i=1,...,k$. If
$(i,i)\in A(H)$, then $V_i$ is a clique, otherwise it is an
independent set.

Suppose, that $V_i = \{v_{l_1},...,v_{l_i}\}$, $l_1 < l_2
<...<l_i$. If $(i,j)\in A(H)$, then $N_{V_j}(v_{l_1})\supseteq
N_{V_j}(v_{l_2})\supseteq...\supseteq N_{V_j}(v_{l_i})$, otherwise
$N_{V_j}(v_{l_i})\supseteq ... \supseteq N_{V_j}(v_{l_1})$.
 So, the partition (\ref{part2}) satisfies the neighborhoods
ordering property.

Let $D$ be a digraph obtained from $H$ by deleting loops and arcs
of the set $\{(i,j): V_i\sim V_j\}$. Then in the digraph $D$

$$N_{out}(i) = \{j : N_{V_j}(v_{l_1})\supseteq
N_{V_j}(v_{l_2})\supseteq...\supseteq
N_{V_j}(v_{l_i})\;and\;neither\; V_i\sim V_j \; nor\; V_i\not\sim
V_j\};$$

$$N_{in}(i) = \{j : N_{V_j}(v_{l_i})\supseteq ... \supseteq N_{V_j}(v_{l_1})\;and\;neither\; V_i\sim V_j \; nor\; V_i\not\sim V_j\}.$$

So, $D$ is a realization of $S(V_1,...,V_k)$.

It remains to show, that $(v_1,...,v_n)$ is the acyclic ordering
of $F = F(V_1,...,V_k)$. All arcs with both ends in $V_l$,
$l=1,...,k$, have the form $(v_i,v_{i+1})$. So, let us consider
$v_i\in V_l$, $v_j\in V_s$, $l\ne s$ such, that $v_i$ and $v_j$
are adjacent in $F$. By the definition of $F$ neither $V_l\sim
V_s$ nor $V_l\not\sim V_s$. Then $l$ and $s$ are adjacent in $H$.
Let $(l,s)\in A(H)$. If $(v_i,v_j)\in A(F)$, then $v_iv_j\in
E(G)$, which could be only if $i < j$. If $(v_j,v_i)\in A(F)$,
then $v_iv_j\not\in E(G)$, which could be only if $j < i$. The
theorem is proved.
\end{pf}

\begin{rem}
If the partition (\ref{part}) is given, it could be tested in a
polynomial time, if it satisfies the conditions of the Theorem
\ref{HThreshcharact}. In case of the positive answer, the proofs
of the Lemma \ref{digraphical} and Theorem \ref{HThreshcharact}
contain the algorithm for reconstruction of the graph $H$ such
that $G$ is $H$-threshold graph.

\end{rem}

The definition of the digraph $F(V_1,...,V_k)$ depends on the
realization $D$ of the family $S(V_1,...,V_k)$. But the family
$S(V_1,...,V_k)$ can have different realizations. The next
proposition shows, that from the point of view of the Theorem
\ref{HThreshcharact} it does not matter, which realization to
choose.

\begin{prop}
{\it Let $D_1$, $D_2$ be two realizations of $S(V_1,...,V_k)$ for
a partition (\ref{part}). If $F_{D_1}(V_1,...,V_k)$  is acyclic,
then $F_{D_2}(V_1,...,V_k)$ is also acyclic.}

\end{prop}

\begin{pf} Suppose, that $F_{D_1}(V_1,...,V_k)$  is acyclic. By the
corollary from the Lemma \ref{digraphical} $D_1$ and $D_2$ have
the same sets of connected components. It follows from the
definition, that $\{i_1,...,i_j\}$ is a connected component of
$D_l$ if and only if $V_{i_1}\cup...\cup V_{i_j}$ is a connected
component of $F_{D_l}(V_1,...,V_k)$, $l=1,2$. So, the definition
of $F$ and the Corollary \ref{observrealiz} imply, that
$F_{D_2}(V_1,...,V_k)$ could be obtained from
$F_{D_l}(V_1,...,V_k)$ by the reversal of all arcs of some of its
connected components. So, $F_{D_2}(V_1,...,V_k)$ is acyclic.
\end{pf}

\section{Graphs with $ThrWidth(G)\leq 2$}
\label{HThresh2}

It is clear, that graphs with $ThrWidth(G) = 1$ are exactly
complete and empty graphs. For every threshold graph $G$
$ThrWidth(G) \leq 2$. But the set of graphs with $ThrWidth(G) \leq
2$ is not reduced to the threshold graphs. For example, on the
figure 2 we can see, that $C_4$ and $P_4$ have the threshold-width
2.

\begin{prop}\label{ThrBipChain}{\it $ThrWidth(G)\leq 2$ if and only if $G$ or $\overline{G}$ is
either threshold, or difference.}

\end{prop}

\begin{pf} By the Theorem \ref{HThreshcharact} the necessity is
straightforward, so let us prove the sufficiency. Let us use the
Theorem \ref{HThreshcharact}.   By the definition there exists the
partition $V(G)=V_1\cup V_2$ such that $V_1$($V_2$) is either
clique or independent set.

This partition satisfies the neighbourhoods ordering property. It
is clear, that the realization of the family $S(V_1, V_2)$ is
either empty digraph (if $V_1\sim V_2$ or $V_1\not\sim V_2$) or
the digraph $D$ with $A(D) = \{(1,2)\}$.

Let us prove that $F = F(V_1,V_2)$ is acyclic. If $V_1\sim V_2$ or
$V_1\not\sim V_2$, then $F$ is empty. Otherwise let $A(D) =
\{(1,2)\}$.

Let $V_1 = \{u_1,...,u_r\}$, $V_2 = \{v_1,...,v_s\}$, where
$N_{V_2}(u_1)\supseteq N_{V_2}(u_2)\supseteq...\supseteq
N_{V_2}(u_r)$, $N_{V_1}(v_s)\supseteq
N_{V_2}(v_{s-1})\supseteq...\supseteq N_{V_2}(v_1)$. Then all arcs
of $F$ with both ends in $V_1$ ($V_2$) have the form
$(u_i,u_{i+1})$, $i=1,...,r-1$ ($(v_i,v_{i+1})$, $i=1,...,s-1$).
Therefore if there exists a directed cycle in $F$, it should
contain arcs $(u_j,v_l)$, $(v_p,u_i)$, $i\leq j$, $l\leq p$ (since
$F$ contains no loops we may assume without lost of generality,
that $i\ne j$). By the definition of $F$, it means that $u_jv_l\in
E(G)$, $u_iv_p\not\in E(G)$. Since $N_{V_2}(u_i)\supseteq
N_{V_2}(u_j)$ we have $u_iv_l\in E(G)$. If $l=p$, then we have the
contradiction. If $l\ne p$ then, as $N_{V_1}(v_p)\supseteq
N_{V_1}(v_l)$, we again have $u_iv_p\in E(G)$. This contradiction
finishes the proof. \end{pf}

\begin{cor}
{\it The class of difference graphs coincides with the class of
$H'$-threshold graphs, where $V(H')=\{1,2\}$, $A(H')=\{(1,2)\}$.}
\end{cor}

\begin{pf}
All $H'$-threshold graphs are difference graphs by the definition
of $\circ_{H'}$ and by Theorem \ref{differgrdef}. Let us show,
that all difference graphs are $H'$-threshold. It is sufficient to
consider connected difference graph $G$ with the bipartition
$(A,B)$ (if $G$ is disconnected, then it is a disjoint union of a
connected difference graph $F$ and $r$ isolated vertices. If $T$
is the threshold representation with respect to $\circ_{H'}$ of
$F$, then $G=(K^2_2\circ_{H'}...\circ_{H'}K^2_2)\circ_{H'}T$ ($r$
multipliers in parentheses)). If $G$ is complete bipartite, then
$G=K_{m,n} = (K_1^2\circ_{H'}...\circ_{H'}
K_1^2)\circ_{H'}(K_2^2\circ_{H'}...\circ_{H'} K_2^2)$ ($m$ and $n$
multipliers in each parentheses). So let further $G$ is not
complete bipartite, which implies, that $|A|,|B|\geq 2$. Then if
$G$ is $H$-threshold, $|H|\leq 2$, then $H$ has no loops. If
$A(H)=\emptyset$, then $G=O_n$, and if $A(H)=\{(1,2),(2,1)\}$,
then $G=K_{m,n} $. So, $H=H'$.
\end{pf}

\begin{figure}[h]
\begin{center}
\includegraphics*[width=125mm] {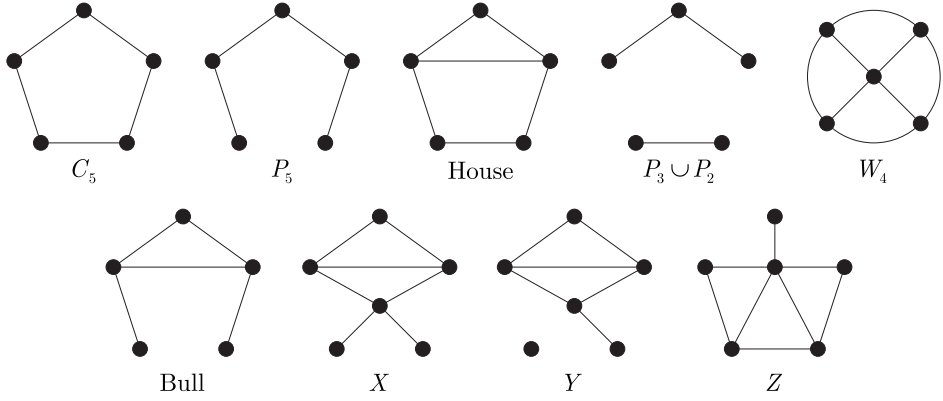}
\caption{The set $L$}%
\end{center}
\label{forbidden}
\end{figure}

\begin{thm}\label{thrdim2}
{\it Let $G$ be a graph. $ThrWidth(G)\leq 2$ if and only if
neither $G$ nor $\overline{G}$ contains one of the graphs from the
set $L = \{C_5,P_5, House, P_3\cup P_2, W_4, Bull, X, Y, Z\}$ as
an induced subgraph.}
\end{thm}

\begin{pf} It is straightforward to check, that every graph from the
set $L$ do not satisfy the Proposition \ref{ThrBipChain}. So we
will prove the sufficiency.

Let us prove firstly, that $G$ is either split, or bipartite, or a
complement of bipartite. After that we will prove, that for each
its part the neighborhoods of its vertices in the another part are
ordered by inclusion.

 Suppose, that neither $G$ nor
$\overline{G}$ is bipartite. We will show, that $G$ is split.

Let $A$ be a maximum clique of $G$ and such, that a subgraph
induced by the set $B=V(G)\setminus A$ have the smallest possible
number of edges. We will prove, that $B$ is an independent set.

Suppose the contrary, i.e. there exist $x,y\in B$ such, that
$x\sim y$. Since $A$ is maximum, there exist vertices of $A$,
which are not adjacent to $x$ ($y$). If all vertices of $A$,
except, possibly, one vertex $u$, adjacent to both $x$ and $y$,
then $A\setminus \{u\}\cup \{x,y\}$ is a clique, which contradicts
the maximality of $A$. So, there exist $u,v\in A$ such, that
$u\not\sim x$, $v\not\sim y$.

It is easy to see, that $|A|\geq 3$. Indeed, if $|A|=2$, then $G$
is triangle-free. It, together with the fact, that $G$ is $\{C_5,
P_5\}$-free, imply that $G$ doesn't contain odd cycles.

Let $w\in A\setminus \{u,v\}$. Because $\overline{G}$ is not
bipartite, there exists $z\in B\setminus \{x,y\}$ such, that
$z\not\sim y$ or $z\not\sim x$. We may assume, that $w\not\sim z$,
since $A$ is a maximum clique.

Let us call the induced cycle $C=C_4$ {\it bad}, if there exists a
vertex $a\in V(G)\setminus C$ such, that $|N(a)\cap C|\geq 2$. By
the assumption of the theorem $G$ does not contain bad $C_4$'s.

If $u\sim y$ and $v\sim x$, then $G$ contains bad $C_4$. Therefore
the following cases are possible: 1) $u\not\sim y$, $v\not\sim x$
and 2) $u\sim y$, $v\not\sim x$. Consider those cases.

1) $u\not\sim y$, $v\not\sim x$.

Let without lost of generality $z\not\sim y$. If $z\sim x$, then
without lost of generality $z\sim v$ (since $G[u,v,y,x,z]\ne
P_3\cup P_2$). As $G[y,x,z,v,w]\ne P_5, C_5$, $w\sim x$. But then
$\{w,v,z,x\}$ form bad $C_4$.

So it is proved, that $z\not\sim x$. Moreover, it is shown, that
for every $t\in B\setminus \{x,y\}$ $t\sim \{x,y\}$ or $t\not\sim
\{x,y\}$.

Let $T_1 = \{t\in B\setminus \{x,y\} : t\sim \{x,y\}\}$, $T_2 =
\{t\in B\setminus \{x,y\} : t\not\sim \{x,y\}\}$. We know from the
considerations above, that $T_2\ne \emptyset$.

Let $t\in T_2$. As $G[u,v,y,x,t]\ne P_3\cup P_2$, without lost of
generality $t\sim v$. Then, since $G[t,v,u,y,x]\ne P_3\cup P_2$,
$t\sim u$. So, we have $T_2\sim \{u,v\}$.

\begin{lem}\label{cl1}
{\it For every $q\in A\setminus \{u,v\}$ $q\sim T_2$ or $q\sim
\{x,y\}$. Moreover, $T_2$ is a clique.}
\end{lem}

\begin{pf} Suppose, that there exists $t\in T_2$ such, that $q\not\sim
t$. The statement, that $q\sim \{x,y\}$ follows from the fact,
that $G[t,v,q,y,x]\ne P_3\cup P_2, P_5$. If there exist
$t_1,t_2\in T_2$ such, that $t_1\not\sim t_2$, then
$G[t_1,v,t_2,y,x]=P_3\cup P_2$. \end{pf}

Let $Q_1=\{q\in A\setminus \{u,v\} : q\sim T_2\}$, $Q_2 =
(A\setminus\{u,v\})\setminus Q_1$. By the Lemma \ref{cl1} $Q_2\sim
\{x,y\}$. Moreover, as $A$ is maximal clique, $Q_2\ne \emptyset$.

\begin{lem}\label{cl2}
{\it $Q_2\sim T_1$. Moreover, $T_1$ is a clique.}
\end{lem}

\begin{pf} Suppose, that there exist $t_1t_2\in T_1$ such, that
$t_1\not\sim t_2$. Since $G[u,v,y,t_1,t_2]\ne P_3\cup P_2$,
without lost of generality $t_2\sim v$. Then either $t_2\sim u$ or
$t_1\sim v$, because $G[u,v,t_2,y,t_1]\ne P_5,C_5$. But
$t_1\not\sim v$, because otherwise $v,t_2,x,t_1$ form bad $C_4$.
So $t_1\not\sim v$, $t_2\sim u$. Analogously, it is easy to see,
that $t_1\not\sim \{u\}$.

By the maximality of the clique $A$, there exists $q\in A\setminus
\{u,v\}$ such, that $q\not\sim t_2$. As $G[q,v,t_2,y,t_1]\ne P_5,
C_5$, $q\sim x$. But then $G[q,v,t_2,x]$ is a bad $C_4$. So it is
proved, that $T_1$ is a clique.

Let us show now, that $T_1\sim Q_2$. Suppose the contrary, i.e.
let there exist $t\in T_1$, $q\in Q_2$ such, that $t\not\sim q$.
By the definition of $Q_2$ there exist $z\in T_2$ such, that
$q\not\sim z$. Since $G[z,u,q,x,t]\ne P_5, C_5$, $t\sim u$. But
then $G[u,q,x,t]$ is a bad $C_4$. \end{pf}

By Lemma \ref{cl1} and Lemma \ref{cl2} $V_1 = Q_2\cup T_1\cup
\{x,y\}$ and $V_2 = Q_1\cup T_2\cup \{u,v\}$ are cliques, $V_1\cup
V_2 = V(G)$. The contradiction with the fact, that $\overline{G}$
is not bipartite, is obtained. So, the case 1) is considered.

2) $u\sim y$, $v\not\sim x$.

\begin{lem}\label{claim3-}
{\it For every $z\in B\setminus\{x,y\}$ $z\sim \{x,y\}$ or
$z\not\sim\{x,y\}$.}
\end{lem}

\begin{pf} Assume, in contrary, that there are exist $z\in B\setminus
\{x,y\}$ such, that the lemma is not satisfied for it.

Let $z\sim x$, $z\not\sim y$. Since $G[v,u,y,x,z]\ne P_5, C_5$,
then $z\sim u$. Consider $w\in A\setminus \{u,v\}$. As
$G[y,x,z,v,w]\ne P_3\cup P_2$, there are edges between $\{v,w\}$
and $\{x,y,z\}$. But it means, that $G[u,y,x,z]$ is a bad $C_4$.

So, $z\sim y$, $z\not\sim x$. Suppose, that $z\sim v$. Then $z\sim
u$ (because $G[z,y,u,v]$ is not a bad $C_4$). Therefore by the
maximality of $A$ there exists $w\in A$ such that $w\not\sim z$.
For this vertex we have $w\sim y$ (as $G[x,y,z,v,w]\ne P_5, C_5$),
and it implies, that $G[w,v,z,y]$ is a bad $C_4$.

So, $z\not\sim v$. But then $z\not\sim u$ (otherwise $G[z,y,u,v,x]
= Bull$). Since $G[x,y,z,w,v]\ne P_3\cup P_2$, there exist some of
the edges from the set $\{wx,wy,wz\}$.

Suppose, that $w\sim x$. Then $w\sim y$ (because otherwise
$G[w,x,y,u]$ is a bad $C_4$). It implies, that $w\sim z$ (as
$G[w,y,x,z,v]\ne Bull$). But then $G[x,y,z,u,v,w]=\overline{Y}$.

Thus $w\not\sim x$. If $w\sim z$, then $w\sim y$ (since
$G[w,u,y,z]$ is not bad $C_4$). It implies, that $G[w,z,y,v,x] =
Bull$.

So, $w\not\sim z$. Then $w\sim y$ and $G[w,u,v,y,z,x] = X$.
\end{pf}

Let $B\setminus \{x,y\} = S_1\cup S_2$, $S_1=\{z\in B : z\sim
\{x,y\}\}$, $S_2=\{z\in B : z\not\sim \{x,y\}\}$. Since
$\overline{G}$ is not bipartite, $S_2\ne \emptyset$.

\begin{lem}\label{claim3}
{\it For every $r\in A\setminus\{u,v\}$ $r\sim \{x,y\}$ or $r\sim
S_2$.}
\end{lem}

\begin{pf} Assume, that there exists $z\in S_2$ such, that $r\not\sim
z$.

Let $z\sim v$. As $G[z,v,r,y,x]\ne P_3\cup P_2$, $r\sim y$ or
$r\sim x$. The situation, when $r\sim x$ and $r\not\sim y$, is
impossible, because otherwise $G[r,u,y,x]$ is a bad $C_4$. If
$r\sim y$, then $r\sim x$ (because $G[x,y,r,v,z]\ne P_5$).

It remains to consider the case, when $z\not\sim v$.  Then $r\sim
y$ or $r\sim x$, since $G[r,v,y,x,z]\ne \overline{W_4}$. As above,
the case, when $r\sim x,$ $r\not\sim y$, is impossible. So $r\sim
y$. As $G[v,u,r,y,x,z]\ne Y$, $z\sim u$ or $r\sim x$. The
situation, when $z\sim u$, $r\not\sim x$ contradicts the fact,
that $G[r,u,y,x,r]\ne Bull$. So $r\sim x$. \end{pf}

Let $A\setminus \{u,v\} = R_1\cup R_2$, $R_1 = \{r\in A\setminus
\{u,v\} : r\sim S_2\}$, $R_2 = (A\setminus \{u,v\})\setminus R_1$.
By the Lemma \ref{claim3} $R_2 \sim \{x,y\}$.

\begin{lem}\label{claim4}
{\it $S_2\sim \{u,v\}$. Moreover, $S_2$ is a clique.}
\end{lem}

\begin{pf} Let us first the first statement of the lemma. Let $z\in
S_2$. Assume, that $z\not\sim v$. We will show, that it is
impossible.

Suppose, that there exists $r\in A\setminus \{u,v\}$ such that
$r\not\sim y$. By the Lemma \ref{claim3} $r\sim z$. Then $r\sim
x$, since $G[z,r,v,y,x]\ne P_3\cup P_2$. But then $G[r,u,y,x]$ is
a bad $C_4$.

So, it is proved that $y\sim A\setminus \{v\}$. Therefore there
exists $s\in B\setminus\{x,y,z\}$ such that $s\sim v$ and
$s\not\sim y$. Indeed, if, on the contrary, $N_B(v)\subseteq
N_B(y)$, then $A' = (A\setminus \{v\})\cup \{y\}$ is a maximum
clique and for the subgraph, induced by the set $B' =
V(G)\setminus A'$, we have $|E(G[B'])| < |E(G[B])|$. It
contradicts the definition of the clique $A$.

As $G[x,y,u,v,s]\ne P_5, C_5$, $s\sim u$. Moreover, $s\not\sim x$,
(because otherwise $G[x,y,u,s]$ is a bad $C_4$) and $s\sim z$
(because otherwise $G[v,s,y,x,z] = \overline{W_4}$). But then
$G[z,s,v,y,x] = P_3\cup P_2$.

So, $z\sim v$. Then $z\sim u$ (see the proof of the Lemma
\ref{claim3}).

Now it is easy to see, that $S_2$ is a clique. Indeed, if there
exist $s_1,s_2\in S_2$ such that $s_1\not\sim s_2$, then
$G[s_1,v,s_2,y,x] = P_3\cup P_2$. \end{pf}

In particular, Lemma \ref{claim4} and the maximality of $A$ imply,
that $R_2\ne \emptyset$.

\begin{lem}\label{claim5}
{\it $R_2\sim S_1$. Moreover, $S_1$ is a clique.}
\end{lem}

\begin{pf} Let there exist $r\in R_2$ and $s\in S_1$ such that
$r\not\sim s$. By Lemma \ref{claim3} $r\sim \{x,y\}$. By the
definition there exists $z\in S_2$ such, that $z\not\sim r$. Lemma
\ref{claim4} implies, that $z\sim \{u,v\}$. Since $G[z,v,r,x,s]\ne
P_5, C_5$, either $z\sim x$ or $s\sim v$. But in the first case
$G[z,v,r,x]$ is a bad $C_4$, and in the second case $G[v,r,x,s]$
is a bad $C_4$. So, it is proved, that $R_2\sim S_1$.

Let us show now, that $S_1$ is a clique. Suppose that there exist
$z_1, z_2\in S_1$ such, that $z_1\not\sim z_2$. As
$G[z_1,x,z_2,u,v]\ne P_3\cup P_2$, there exists at least one edge
between $\{z_1,z_2\}$ and $\{u,v\}$. At the same time, if $z_1\sim
v$ and $z_1\not\sim u$, then $G[z_1,v,u,y]$ is a bad $C_4$.

So, without lost of generality $z_1\sim u$. Then $z_2\not\sim u$
(because otherwise $G[z_1,x,z_2,u]$ is a bad $C_4$). Since
$G[v,u,z_1,x,z_2]\ne P_5, C_5$, $z_1\sim v$. It implies, that
$z_2\not\sim v$ (otherwise $G[v,z_1,x,z_2]$ is a bad $C_4$).

The maximality of $A$ implies the existence of $w\in A$ such, that
$w\not\sim z_1$. $w\not\sim x$, as $G[w,v,z_1,x]$ is not a bad
$C_4$. But then $G[w,v,z_1,x,z_2] =$ $P_5$ or $C_5$. \end{pf}

By Lemma \ref{claim4} and Lemma \ref{claim5} $V_1 = R_2\cup
S_1\cup \{x,y\}$ and $V_2 = R_1\cup S_2\cup \{u,v\}$ are cliques,
$V_1\cup V_2 = V(G)$. The contradiction with the fact, that
$\overline{G}$ is not bipartite, is obtained. The case 2) is
considered.

So, it is proved, that $G$ or $\overline{G}$ is either split or
bipartite. Let $(A,B)$ be the bipartition of $G$. Let us show,
that the neighborhoods of vertices from $A(B)$ are ordered by
inclusion.

Let us suppose the contrary, i.e. there exist $u,v\in A$, $x,y\in
B$ such, that $u\sim x$, $v\not\sim x$, $u\sim y$, $v\not\sim y$.

Suppose, that $G$ is bipartite. If $|V(G)|=4$, then $\overline{G}$
the statement of the theorem obviously holds. Let there exists
$z\in B\setminus \{x,y\}$. Since $G[u,v,x,y,z]\ne \overline{W_4},
P_3\cup P_2$, $z\sim u,v$. But then $G[u,v,x,y,z]=P_5$. This
contradiction proves the theorem for bipartite graphs.

Taking into account Observation \ref{obs1}, it remains to consider
the case, when $G$ is split and neither bipartite nor a complement
of bipartite.

The following statements hold:

a) $N(x)\cup N(y) = A$ (since $G$ does not contain $Bull$);

b) for every $z\in B\setminus \{x,y\}$ $|N(z)\cap \{u,v\}|\leq 1$
(by the same reason as in a));

c) $|A|\geq 3$, $|B|\geq 3$ (otherwise either $G$ or
$\overline{G}$ is bipartite).

Let $z\in B\setminus \{x,y\}$, $w\in A\setminus \{u,v\}$, $w\sim
x$. As $G[u,v,x,y,w,z]\ne Y,\overline{Z}$, at least one of the
edges $zu,zv,zw$ belongs to $E(G)$. If there exists exactly one of
this edges, then $G[u,v,w,z,y] = Bull$, $G[u,v,w,x,y,z] = X$,
$G[u,v,w,z,y] = Bull$, respectively. Therefore, taking into
account b), either $zw,zv\in E(G)$, $zu\not\in E(G)$ or $zw,zu\in
E(G)$, $zv\not\in E(G)$.

In the first case $w\sim y$ (since $G[w,v,z,y,x] \ne Bull$), which
implies, that $F = G[u,v,x,y,w,z] = \overline{Y}$. In the second
case $w\sim y$ (since $F \ne \overline{X}$), which implies, that
$F =\overline{Y}$.

The theorem is proved \end{pf}

\end{document}